\documentclass[11pt]{article}
\usepackage{amsthm, amsmath, amssymb, amsfonts, url, booktabs, tikz, setspace, fancyhdr, bm}
\usepackage{hyperref}
\usepackage{amsthm}
\usepackage{geometry}
\usepackage{hyperref, enumerate}
\usepackage[shortlabels]{enumitem}
\usepackage[babel]{microtype}
\usepackage[english]{babel}
\usepackage[capitalise]{cleveref}
\usepackage{comment}
\usepackage{bbm}
\usepackage{csquotes}
\usepackage{mathabx}
\usepackage{tikz}
\usepackage{graphicx}
\usepackage{float}
\usepackage{xcolor}
\usetikzlibrary{positioning, arrows.meta, shapes.geometric}
\counterwithin{figure}{section}

\newtheorem{theorem}{Theorem}[section]

\newtheorem{conj}[theorem]{Conjecture}

\theoremstyle{definition}

\newtheorem*{defn-non}{Definition}

\newcommand{\Yemph}[1]{\textcolor{black}{\emph{#1}}}

\newlist{Case}{enumerate}{2}
\setlist[Case, 1]{%
    label           =   {\bfseries Case \arabic*.},
    labelindent=1em ,labelwidth=1.3cm, labelsep*=1em, leftmargin =!
}
\setlist[Case, 2]{%
    label           =   {\bfseries Subcase \arabic{Casei}.\arabic*.},
    labelindent=-1em ,labelwidth=1.3cm, labelsep*=1em, leftmargin =!
}

\newcommand{\C}[1]{{\protect\mathcal{#1}}}

\usepackage{todonotes} 

\newcommand{\eps}{\varepsilon}

\title{Piercing independent sets in graphs without large induced matching} 
\author{
Jiangdong Ai\thanks{School of Mathematical Sciences and LPMC, Nankai University, Tianjin, P.R. China. Email: {jd@nankai.edu.cn}. Supported by NNSFC No.12161141006.}
\and
Hong Liu\thanks{Extremal Combinatorics and Probability Group (ECOPRO), Institute for Basic Science (IBS), Daejeon, South Korea. Emails: {\texttt \{hongliu, zixiangxu\}@ibs.re.kr}. Supported by the Institute for Basic Science (IBS-R029-C4).}
\and
Zixiang Xu\footnotemark[2]
\and
Qiang Zhou\thanks{Academy of Mathematics and Systems Science (AMSS), Chinese Academy of Sciences (CAS), Beijing, P.R. China, and University of Chinese Academy of Sciences (UCAS), Beijing, P.R. China. Email: {zhouqiang2021@amss.ac.cn}. Supported by National Key R\&D Program of China, No.2023YFA1009602, and IBS-R029-C4.}
}
\begin{document}
\date{}
\maketitle
\begin{abstract}
Given a graph $G$, denote by $h(G)$ the smallest size of a subset of $V(G)$ which intersects every maximum independent set of $G$. We prove that any graph $G$ without induced matching of size $t$ satisfies $h(G)\le \omega(G)^{3t-3+o(1)}$. This resolves a conjecture of Hajebi, Li and Spirkl~(Hitting all maximum stable sets in $P_{5}$-free graphs. J. Combin. Theory Ser. B, 165:142–163, 2024).
\end{abstract}

\maketitle

\section{Introduction}
For a graph $G$, let $h(G)$ denote the smallest size of a subset of $V(G)$ intersecting every maximum independent set of $G$. Bollob\'as, Erd\H{o}s and Tuza~\cite{1999BookErdos,1991ErdosCollection} conjectured that any graph $G$ with linear (in $|G|$) independence number must have sublinear $h(G)$. This easy-to-state conjecture is surprisingly difficult and remains wide open. It is then natural to consider special families of graphs. For example, Alon~\cite{2022Alon} suggested to study 3-colorable graphs. Motivated by this and some similarities to another well-studied problem of $\chi$-boundedness, Hajebi, Li and Spirkl~\cite{2023P5Free} very recently consider the problem of bounding $h(G)$ by the clique number $\omega(G)$. Indeed, they observe that (1) given $G$, $h(H)\le \omega(H)$ for every induced $H\subseteq G$ if and only if $G$ is perfect; (2) there are graphs $G$ with arbitrarily large $h(G)$ and girth; and (3) if for any graph $G$ in a hereditary family $\C G$ (i.e.~closed under taking induced subgraphs),  $h(G)$ can be bounded by a polynomial function of $\omega(G)$, then $\C G$ satisfies Erd\H{o}s-Hajnal conjecture~\cite{1989DAMEH}, yet another central conjecture in extremal and structural graph theory, that every graph in $\C G$ has polynomial size homogeneous set. We refer the interested readers to the survey of 
Scott and Seymour~\cite{2018SurveyChiBound} for more about $\chi$-boundedness and the recent work~\cite{2024FoxEquivalent,2023P5ErdosHajnal,2023BVCErdosHajnal} for more on Erd\H{o}s-Hajnal conjecture. 

Hajebi, Li and Spirkl~\cite{2023P5Free} proved that for graphs  $G$ without induced $P_5$, $h(G)\le f(\omega(G))$ for some function $f$. They raised the following conjecture that a similar phenomenon holds for graphs without large induced matching and furthermore one can take a polynomial binding function $f$. An \Yemph{induced matching} of size $t$ consists of vertex set $\{a_{i},b_{i}\}_{i\in [t]}$ and edge set $\{a_{i}b_{i}\}_{i\in [t]}$.

\begin{conj}[\cite{2023P5Free}]\label{conj:Matching}
   Let $G$ be a graph without induced matching of size $t$, then $h(G)\le \omega(G)^{O_t(1)}$. 
\end{conj}

Hajebi, Li and Spirkl~\cite{2023P5Free} resolved the first case $t=2$. Quoting them, `\emph{it is not even known (strikingly enough)}' that whether $h(G)$ can be bounded by any function of $\omega(G)$ for $t=3$.

In this short note, we prove~\cref{conj:Matching}. Our proof is inspired by recent work~\cite{2024GraphToGeom}.

\begin{theorem}\label{theorem1}
   Let $G$ be a graph without induced matching of size $t$, then $h(G)\le 10t^t\omega(G)^{3t-3}\log\omega(G)$.
\end{theorem}

\section{The proof}
For a set system $\mathcal{F}\subseteq 2^{V}$, the \Yemph{transversal number} of $\mathcal{F}$, denoted by $\tau(\C F)$, is the minimum size of subset $T\subseteq V$ such that $T\cap F\neq\emptyset$ for every set $F\in \mathcal{F}$. The \Yemph{fractional transversal number} of $\mathcal{F}$, denoted by $\tau^*(\C F)$, is the minimum of $\sum_{v\in V}g(v)$, taken over all functions $g$: $V\to [0,1]$ such that $\sum_{v\in F}g(v)\geq 1$ for every $F\in \mathcal{F}$. The \Yemph{Vapnik–Chervonenkis dimension} (VC-dimension for short) of $\mathcal{F}$ is the maximum cardinality of a subset $S\subseteq V$ such that for every $S'\subseteq S$ there exists $F\in\mathcal{F}$ with $F\cap S=S'$. Note that by definition, $\tau^*(\C F)\le \tau(\C F)$. The well-known $\eps$-net theorem of Haussler and Welzl~\cite{1987DCG} (also see the book of Matou\v{s}ek~\cite{2002MatousekLecture} and~\cite{1994DingSeymour}) provides an inverse relation for set systems with bounded VC-dimension. 

\begin{theorem}[\cite{1987DCG}]\label{thm:transversal}
Let $\mathcal{F}$ be a set system with VC-dimension $d$, then $\tau(\mathcal{F})\leq 2d\tau^{*}(\mathcal{F})\log(11\tau^{*}(\mathcal{F}))$.
\end{theorem}

Now, let $G$ be an $n$-vertex graph with no induced matching of size $t$. Consider the set system $\C F\subseteq 2^{V(G)}$ consisting of all maximum independent sets of $G$. By definition, $h(G)=\tau(\mathcal{F})$. Considering the constant function $g\equiv 1/\alpha(G)$ over $V(G)$, we see that 
$$\tau^{*}(\C F)\leq n/\alpha(G)\leq\chi(G)\leq \omega(G)^{2t-2},$$ 
where the last inequality is a classical result of Wagon~\cite{1980JCTB} on $\chi$-boundedness for graphs without large induced matching. It remains to bound the VC-dimension of $\C F$.

Suppose the VC-dimension of $\C F$ is $d$ and $S=\{v_{1},v_{2},\ldots,v_{d}\}$ is a set such that for any subset $S'\subseteq S$, there is a maximum independent set $I_{S'}\in \C F$ with $I_{S'}\cap S=S'$. In particular, the set $S$ itself is an independent set by considering $I_S$. For each $i\in [d]$, as $v_i\notin I_{S\setminus \{v_i\}}$, and $I_{S\setminus \{v_i\}}$ is a maximum independent set, there exists some vertex $u_{i}\in I_{S\setminus \{v_i\}}$ such that $v_{i}u_{i}\in E(G)$ and $v_{j}u_{i}\notin E(G)$ for any $j\neq i$. Since $S$ is an independent set, we must have $u_i\notin S$. Moreover, $u_{1},\ldots,u_{d}$ are all distinct. Indeed, if $u_j=u_i$ for some $j\neq i$, then $v_ju_j=v_ju_i\notin E(G)$, a contradiction. Thus, $\{v_{i},u_{i}\}_{i\in [d]}$ forms a matching of size $d$, and $v_{1},v_{2},\ldots,v_{d}$ form an independent set. Then, as $G$ has no induced matching of size $t$, $G[\{u_1,\ldots, u_d\}]$ has no independent set of size $t$, and so by off-diagonal Ramsey~\cite{1935ES}, $d<R(\omega(G)+1,t)\le \binom{\omega(G)+t-1}{t-1}$. Thus by~\cref{thm:transversal}, we have 
   $$h(G) \leq 2\binom{\omega(G)+t-1}{t-1}\omega(G)^{2t-2}\log(11\omega(G)^{2t-2})\le 10t^t\omega(G)^{3t-3}\log\omega(G),$$
completing the proof.

\bibliographystyle{abbrv}
\bibliography{Matching}
\end{document}